\documentclass[]{interact}

\usepackage{epstopdf}

\usepackage{natbib}
\bibpunct[, ]{(}{)}{;}{a}{}{,}

\theoremstyle{plain}

\theoremstyle{definition}

\theoremstyle{remark}

 \usepackage{enumerate}
\usepackage{amsmath,amsfonts,amssymb}
\usepackage{chemformula}
\usepackage{float}

\usepackage{subfig}
\usepackage{color}
\usepackage{bm}
\usepackage{bbm}
\usepackage{lscape}

\usepackage{graphicx}%
\usepackage{multirow}%
\usepackage{amsmath,amssymb,amsfonts}%
\usepackage{amsthm}%
\usepackage{mathrsfs}%
\usepackage{xcolor}%
\usepackage{textcomp}%
\usepackage{manyfoot}%
\usepackage{booktabs}%
\usepackage{algorithm}%
\usepackage[noend]{algorithmic} 

\usepackage{listings}%
\usepackage{array, float, indentfirst, longtable, tikz, xcolor}
\usepackage{caption}
\usepackage{subfig}
\usepackage{tikz}
\usetikzlibrary{patterns}

\allowdisplaybreaks

\newcommand{\revise}[1]{{\color{black}#1}}

\begin{document}

\articletype{ARTICLE TEMPLATE}

\title{Discretization-Based Solution Approaches  for the Circle Packing Problem\footnote{This research is based on the first author's master thesis \citep{tacspinar2021discretization}, and its four-page summary has appeared in the Proceedings of the 17th International Symposium on Operations Research in Slovenia.}}

\author{
\name{Rabia Ta{\c{s}}p{\i}nar\textsuperscript{a} \and Burak Kocuk\textsuperscript{a}\thanks{CONTACT Burak Kocuk. Email: burakkocuk@sabanciuniv.edu}}
\affil{\textsuperscript{a}Industrial Engineering Program,  Sabanc{\i} University, Istanbul, Turkey 34956}
}

\maketitle

\begin{abstract}
The problem of packing a set of circles into the smallest surrounding container is considered. This problem arises in different application areas such as automobile, textile, food, and chemical industries. The so-called circle packing problem can be cast as a nonconvex quadratically constrained program, and is difficult to solve in general. An iterative solution approach based on a bisection-type algorithm on the radius of the larger circle is provided. The present algorithm discretizes the container into small cells and solves two different integer linear programming formulations proposed for a restricted and a relaxed version of the original problem. The present algorithm is enhanced with solution space reduction, bound tightening and variable elimination techniques. Then, a computational study is performed to evaluate the performance of the algorithm. The present algorithm is compared with BARON and Gurobi that solve the original nonlinear formulation and heuristic methods from literature, and obtain promising results. 
\end{abstract}


\section{Introduction}
\label{sec:intro}

The circle packing problem (CPP) is concerned with packing a given number of different circles into a larger container, such as a square, a rectangle, or a circle, in such a way that circles do not overlap and each circle is entirely in the container. This well-known NP-Hard    problem arises in different application areas including packing circular shaped objects  into the smallest box \citep{LiteratureReview_Castillo2008}, some applications from nanotechnology, telecommunication, electrical, oil, automobile industries \citep{Application_fibreoptic_Huang2002, LiSun_chemistry2009, MinRadiusofContainer_Sugihara2004DiskPF}, forestry   \citep{LiteratureReview_Hifi2009}, location analysis \citep{Application_facilitylayout_Castillo2010}  \revise{and social distancing \citep{bortolete2022support}}. 
\revise{Recent literature also focuses on the related problems including balanced circle packing problem \citep{romanova2022balanced},  proportional circle packing problem \citep{romanova2023proportional} and  circle bin packing problem \citep{yuan2022adaptive,tole2023simulated}.}

Different objectives are considered depending on the application setting in the literature, such as minimizing the area of the surrounding container, maximizing the number of circles packed into a fixed-size container \citep{MaxNumber_Galiev2013}, or maximizing the minimum distance between any two circles \citep{MaxMinDist_maxradii7_Locatelli2002}.
In this study, CPP with the objective of minimizing the radius of the surrounding circle is primarily considered although the present approach can be applied to other containers as well. 

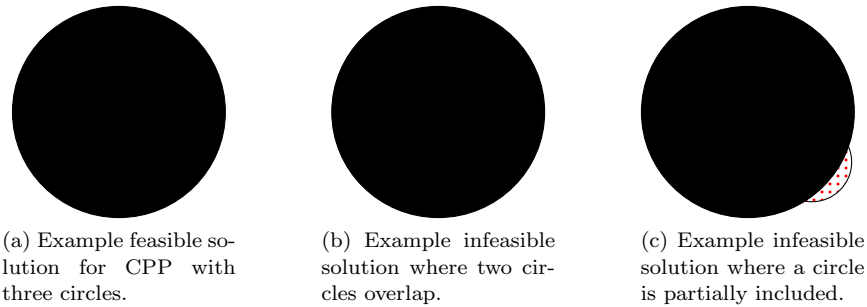
\begin{figure}[H]
\centering
\subfloat[Example feasible solution for CPP with three circles.]{\label{FirstExample_picture_parta}
\centering
\begin{tikzpicture}[scale=0.7]


\filldraw[pattern=north east lines] (2.15,2.85) circle (1cm);

\filldraw[pattern color=blue, pattern=crosshatch] (2,1.2) circle (0.5cm);

\filldraw[pattern color=red, pattern=dots] (3.5,1.5) circle (0.75cm);

\filldraw[fill opacity=0.1] (2.5,2.15) circle (2cm);
\end{tikzpicture}
}
\hspace{1cm}
\subfloat[Example infeasible solution where two circles overlap.]{\label{FirstExample_picture_partb}
\centering\begin{tikzpicture}[scale=0.7]


\filldraw[pattern=north east lines] (2.75,2.75) circle (1cm);

\filldraw[pattern color=blue, pattern=crosshatch] (2,1.2) circle (0.5cm);

\filldraw[pattern color=red, pattern=dots] (3.5,1.5) circle (0.75cm);

\filldraw[fill opacity=0.1] (2.5,2.15) circle (2cm);
\end{tikzpicture}
}
\hspace{1cm}
\subfloat[Example infeasible solution where a circle is partially included.]{\label{FirstExample_picture_partc}\centering\begin{tikzpicture}[scale=0.7]


\filldraw[pattern=north east lines] (2.15,2.85) circle (1cm);

\filldraw[pattern color=blue, pattern=crosshatch] (2,1.2) circle (0.5cm);

\filldraw[pattern color=red, pattern=dots] (3.7,1.2) circle (0.75cm);

\filldraw[fill opacity=0.1] (2.5,2.15) circle (2cm);
\end{tikzpicture}
}
\caption{Example feasible and infeasible placements of circles for CPP.} 
\label{SolutionExamples_picture}

\end{figure}

Here is the precise mathematical formulation of CPP: Denote the set of circles by $\mathcal{C}$, and the radius of circle $c\in \mathcal{C}$  by $r_c$. Assume that the center of the surrounding circle is located at the origin, and $R$ is a decision variable denoting its radius. The center of circle $c$ is represented by the decision variables $(x_c,y_c)$. Then, CPP can be modelled as  the following nonconvex quadratically constrained program:
\begin{subequations}\label{quadForm}
  \begin{align}
 \min \  & R  \\
 \text{s.t.}  \  & {(x_c - x_k)^2 + (y_c - y_k)^2} \ge (r_c + r_k)^2 \qquad c,k\in \mathcal{C} :\; c\neq k \label{quadconst1} \\
 \   & {x_c^2 + y_c^2 } \le (R - r_c)^2 \qquad\qquad\qquad\qquad\ \;\,\, c\in \mathcal{C} \label{quadconst2}\\
  \  &x_c \in \mathbb{R}, \, y_c\in \mathbb{R}, \, R \ge r_c  \qquad \qquad\qquad\quad\;\; c\in \mathcal{C} . \label{quadconst3}
  \end{align}
\end{subequations}
Under the stated feasibility rules of the problem, Figure \ref{FirstExample_picture_parta} shows a feasible suboptimal solution for an instance of the problem consisting of three circles with radii 0.5, 1, 0.75 and 0.5 units where the surrounding container's radius is 2 units. Two possible infeasible configurations are also given: Figure \ref{FirstExample_picture_partb} includes two overlapping circles violating constraint \eqref{quadconst1}; and Figure \ref{FirstExample_picture_partc} is an example where a circle is not fully contained in the container, violating constraint \eqref{quadconst2}.

The vast majority of the CPP literature focuses on developing heuristic methods, which involve constructing a feasible packing of the circles, and then improving it with a search algorithm, see e.g., \cite{minRadius_Huang2006, minRadius_Francesco2014, minRadius_Zeng2016, torres2020binary, LIU2024}. Another stream of research focuses on nonlinear programming techniques to find high quality packings, see e.g., \cite{GlobalSolver_Stoyan2008, GlobalSolver_HUANG2011, GlobalSolver_KunHe_2015, lai2022iterated, ikebe2023mixed}. 
However, none of the studies certify the optimality of the solutions obtained, which are only compared against the best-known results from the literature. 
Hence, the use of systematic global optimization approaches to solve CPP is lacking in the literature. The reason is that it is quite challenging to solve the nonlinear formulation~\eqref{quadForm} with global optimization solvers directly, as confirmed by the preliminary experiments. 

In this article, an iterative solution approach  based on a bisection-type algorithm considering the minimum radius objective is proposed, which  converges to a global optimal solution of CPP. The present solution procedure relies on discretizing the container into smaller cells, and iteratively solves two integer linear programming (ILP) formulations designed for restricted and relaxed versions of the original problem. This allows utilizing mature ILP solvers in order to certify lower and upper bounds for the optimal value of problem~\eqref{quadForm} efficiently. The closest works in the literature to the present article are \cite{MaxNumber_Litvinchev2014, MaxNumber_Litvinchev2015_CircularLike}, which also use discretization but  primarily work with LP relaxations enhanced with valid inequalities as opposed to the present ILP-based approach.

{In the present approach, ILP formulations are constructed based on the discretization of the surrounding circle into smaller squares. The corner points of squares are the guiding points within the ILP formulations, and these points can be represented by utilizing a logarithmic number of binary decision variables in the grid size.  Apart from the two ILP formulations developed for the restricted  and relaxed versions of CPP,  some sub-methods are proposed to start with a better lower bound within the proposed algorithm. Then,   an algorithm based on continuous solution space is developed to shrink the solution space for locating the center of each circle   iteratively. By using this algorithm,  a portion of decision variables are eliminated without any further investigation by pre-processing. These enhancements significantly improve the performance of the present algorithm and \revise{help  solve instances consistently faster than global solvers such as BARON and Gurobi}. \revise{They also enable   obtaining comparable or better feasible solutions compared to heuristic methods from the literature along with optimality guarantees}.
}

{The remainder of the article is organized as follows.  In Section \ref{sec:ILPformulations},  ILP formulations developed for the restricted and relaxed versions of the original problem~\eqref{quadForm} are introduced. In Section \ref{sec:method},  the proposed discretization-based solution approach  is presented as well as the algorithmic enhancements designed for this algorithm. Then, Section \ref{sec:computational} contains a computational study analyzing the performance of the present algorithm and the proposed ILP formulations. Finally, the article is concluded in Section \ref{sec:conclusion}.
}

\section{Integer Linear Programming Formulations}\label{sec:ILPformulations}

In this section,  the ILP formulations developed for  CPP are presented. It is assumed that the surrounding container is circular-shaped; however, the present approach can be applied to other containers, \revise{which can be outer-approximated by a rectangle,} as well.  

Recall the nonconvex quadratically constrained program~\eqref{quadForm}. 
In an optimal solution of the this formulation, the centers of circles can be located into anywhere in $\mathbb{R}^2$ space. However, this problem is hard to solve in the continuous space by using the global solvers BARON and Gurobi. Hence, the continuous solution space is discretized into smaller squares to design the algorithm. This discretization scheme is accompanied by a restricted and a relaxed version under the assumption that a candidate radius for the surrounding circle is given as~$R$. In the restricted version,  the corners of each cell are considered as candidate points to locate the centers of circles. If feasible, this solution will give a feasible solution for CPP, and an upper bound of at most $R$ for problem~\eqref{quadForm}. In the relaxed version,  the centers of circles to be located are allowed to be at any point in a cell. In this version, any pair of circles are allowed to overlap in a well-defined and limited manner so that the resulting model gives a systematic way to construct a relaxation for CPP. If infeasible, the relaxed model will certify a lower bound of~$R$ for problem~\eqref{quadForm}. The details of each formulation are given in the following subsections.

\subsection{Restricted Version}\label{rest_logver}

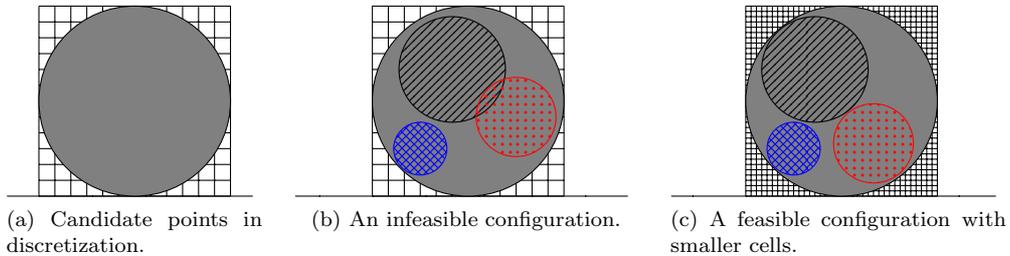
\begin{figure}[H]
\centering

\subfloat[Candidate points in discretization.]{\label{cell_restriction_example_picture1}\centering\begin{tikzpicture}[scale=0.7]
\draw[step=1.0, draw opacity = 0.0] (-0.6,0) grid (4.0,0);

\draw[step=0.3,black,thin,draw opacity=0.15] (0,0) grid (3.6,3.6);
    \foreach \x in {0.9,1.2,1.5,1.8,2.1,2.4,2.7}
    {
        \foreach \y in {0.9,1.2,1.5,1.8,2.1,2.4,2.7}
        {
            \node[draw, circle, inner sep=0.015cm, fill] at (\x,\y) {};
        }
    }
    \foreach \x in {0.6,3.0}
    {
        \foreach \y in {1.2,1.5,1.8,2.1,2.4}
        {
            \node[draw, circle, inner sep=0.015cm, fill] at (\x,\y) {};
        }
    }
    \foreach \x in {1.2,1.5,1.8,2.1,2.4}
    {
        \foreach \y in {0.6,3.0}
        {
            \node[draw, circle, inner sep=0.015cm, fill] at (\x,\y) {};
        }
    }
\draw[fill=gray, opacity = 0.15] (1.8,1.8) circle (1.8 cm);
\end{tikzpicture}
}\hspace{0.3cm}
\subfloat[An infeasible configuration.]{\label{cell_restriction_example_picture2}\centering\begin{tikzpicture}[scale=0.7]
\draw[step=1.0, draw opacity = 0.0] (-1.45,0) grid (4.8,0);

\draw[step=0.3,black,thin,draw opacity=0.15] (0,0) grid (3.6,3.6);

\draw[fill=gray, opacity = 0.15] (1.8,1.8) circle (1.8 cm);

\filldraw[black, pattern color=black, opacity=0.6, pattern=north east lines, draw opacity = 1, fill opacity = 1] (1.5,2.4) circle (1.0cm);
\filldraw[red, pattern color=red, opacity=0.6, pattern= dots, draw opacity = 1, fill opacity = 1] (2.7,1.5) circle (0.75cm);
\filldraw[blue, pattern color=blue, opacity=0.6, pattern=crosshatch, draw opacity = 1, fill opacity = 1] (0.9,0.9) circle (0.5cm);
\end{tikzpicture}
}\hspace{0.3cm}
\subfloat[A feasible configuration with smaller cells.]{\label{cell_restriction_example_picture3}\centering\begin{tikzpicture}[scale=0.7]
\draw[step=1.0,draw opacity = 0.0] (-1.4,0) grid (4.7,0);

\draw[step=0.1,black,thin,draw opacity=0.15] (0,0) grid (3.6,3.6);

\draw[fill=gray, opacity = 0.15] (1.8,1.8) circle (1.8 cm);

\filldraw[black, pattern color=black, opacity=0.6, pattern=north east lines, draw opacity = 1, fill opacity = 1] (1.3,2.4) circle (1.0cm);
\filldraw[red, pattern color=red, opacity=0.6, pattern=dots, draw opacity = 1, fill opacity = 1] (2.4,1.0) circle (0.75cm);
\filldraw[blue, pattern color=blue, opacity=0.6, pattern=crosshatch, draw opacity = 1, fill opacity = 1] (0.9,0.9) circle (0.5cm);
\end{tikzpicture}
}
\caption{Example discretizations of the circle and example placements of circles.}
\label{RestrictedDiscreatizationExample_picture}
\end{figure}

The main features of the formulation can be motivated with the running example from Figure~\ref{SolutionExamples_picture}. First, divide the smallest square containing the candidate surrounding circle with radius~$R=1.8$ into  square-shaped cells as in Figure~\ref{cell_restriction_example_picture1}.  By using the candidate points, subsets of the bullets in Figure \ref{cell_restriction_example_picture1}, 
it is possible to come up with a restricted formulation, which models all constraints of the original problem \revise{with the following property: If the restriction has a feasible solution, then this solution is also feasible for the original problem}. 
In Figures \ref{cell_restriction_example_picture2} and \ref{cell_restriction_example_picture3}, the circles are tried to be packed into a surrounding circle with radius $R=1.8$ units. The cell sizes in these figures are 0.3 and 0.1 units, respectively. Depending on the granularity of the discretization, it may not be possible to find any feasible configuration for a given radius for the surrounding circle, as in  Figure \ref{cell_restriction_example_picture2}. On the other hand, a finer grid with smaller cells \revise{can} enable  finding a feasible solution for the original problem, 
as exemplified in Figure~\ref{cell_restriction_example_picture3}.

This modeling framework is  formalized as follows: Divide the smallest square containing the candidate surrounding circle with radius~$R$ into small square-shaped cells whose side length is $\delta$ where   $\theta\delta = R$ for some positive integer $\theta$. The corners of these cells are ``guiding points" used in these constructions. It is assumed that the diagonal size of each cell is smaller than the radius of the smallest circle, i.e., $\delta\sqrt{2}<\min\{r_1,\cdots,r_n\}$. If a corner of any cell is out of the surrounding circle,   this point is eliminated from the set of candidate points since the size of each cell is smaller than each circle.

Define the set $\mathcal{I}$ to represent the \revise{$x$-coordinates (respectively, $y$-coordinates) of the corners of the cells,} 
where $|\mathcal{I}|=2\theta+1$ such that $\mathcal{I}=\{0,1,\cdots,2\theta\}$. 
The surrounding circle's center is assumed to be at the origin of the coordinate system, which is represented by $(\theta,\theta)$. Hence, the coordinates of the candidate point represented by $(i,j)$ are given with $\left((i-\theta)\delta,(j-\theta)\delta\right)$.  Then, there is a total of $(2\theta+1)^2$ guiding points within the surrounding circle with radius~$R$, and a subset of these points are candidate points to locate the center of circle $c \in \mathcal{C}$.

\revise{A direct approach to represent candidate points would be to define a pair of binary variables for each such point. In this case, the number of binary variables grows linearly in the grid size, which does not scale well for large problem instances (see, \citep{tacspinar2021discretization}, for preliminary experiments). Instead,} 
 a logarithmic number of binary decision variables in the grid size are defined to locate the center of each circle, which is an idea inspired from \cite{logarithmicVielma, logarithmicshabbir}. For this purpose,    a new parameter  $\Theta:=\lceil\log_2{(2\theta)}\rceil$ is introduced and  two vectors of binary decision variables are defined: $\alpha_{c}\in\{0,1\}^\Theta$ and $\beta_{c}\in\{0,1\}^\Theta$. These binary column vectors correspond to the corner point of cell on $x$-axis ($y$-axis) in base 2, where the center of circle $c$ is located. 
Another decision variable $\psi_{i,j,c}$ is also defined to ensure that circle $c$ is totally included in the containing circle by stating that the circle's center should be located from the set $ \mathcal{L}_c := \{(i,j)\in\mathcal{I}^2: ((\theta-i)\delta)^2+((\theta-j)\delta)^2 \le (R - r_c)^2\}$.

To ensure the feasibility for the original problem, overlapping of two circles should be banned. Hence, the candidate points to place the centers of distinct circles $c$ and $k$  should be far enough. In order to guarantee this, the pairs of the minimum required number of cells on $x$ and $y$ axes between the centers of the circles $c$ and $k$ are considered to avoid overlapping of these circles. For this purpose,  the following set is defined:
\parbox{15cm}{\parindent0.5cm $\mathcal{N}_{c,k}=\Big\{(u_1,u_2): u_1, u_2 \in \mathcal{I},\, u_1,u_2 \ge 0,\, (u_1\,\delta)^2+(u_2\,\delta)^2 \ge (r_c+r_k)^2$, 
\par \parindent2.0cm $\big(u_1\,\delta-\delta\big)^2+(u_2\,\delta)^2 < (r_c+r_k)^2,\, (u_1\,\delta)^2+\big(u_2\,\delta-\delta\big)^2 < (r_c+r_k)^2\Big\}$.}

Before introducing the corresponding formulation,   one more binary variable is needed: $\pi_{u_1,u_2,c,k}$. With the help of $\pi_{u_1,u_2,c,k}$ and the set $\mathcal{N}_{c,k}$, it is possible to ensure that any two circles are non-overlapping by stating the sum of all such variables is equal to 1 for each pair of distinct circles $c$ and $k$, and the decision variable $\pi_{u_1,u_2,c,k}$ can be assigned one if the centers of circles $c, k$ are located at two points which have at least $u_1$ cells on $x$-axis, and $u_2$ cells on $y$-axis between them, i.e., $(u_1,u_2)\in\mathcal{N}_{c,k}$. Finally,  the following formulation is obtained:
\begin{subequations}\label{logrestform}
  \begin{alignat}{3}
 &\sum\limits_{t=1}^{\Theta}{2^{(t-1)}\;(\alpha_{c,t}-\alpha_{k,t})} \ge  \sum_{(u_1,u_2)\in \mathcal{N}_{c,k}} u_1\; \pi_{u_1,u_2,c,k} &\quad& c,k\in \mathcal{C}:c<k\label{logvarmodel_cons1a}\\ 
 &\sum\limits_{v=1}^{\Theta}{2^{(v-1)}\;(\beta_{c,v}-\beta_{k,v})} \ge  \sum_{(u_1,u_2)\in \mathcal{N}_{c,k}} u_2\; \pi_{u_1,u_2,c,k} &\quad& c,k\in \mathcal{C}: c<k\label{logvarmodel_cons1b}\\ 
 &\sum_{(u_1,u_2)\in \mathcal{N}_{c,k}} \pi_{u_1,u_2,c,k} = 1 &\quad& c,k\in \mathcal{C}: c<k\label{logvarmodel_cons1c}\\ 
 & \sum\limits_{t=1}^{\Theta}{2^{(t-1)}\;\alpha_{c,t}} \le \sum\limits_{i\in\mathcal{I}:\exists j, (i,j)\in\mathcal{L}_c}{i \; \psi_{i,j,c}} &\quad& c\in \mathcal{C} \label{logvarmodel_cons2a}\\  
 & \sum\limits_{v=1}^{\Theta}{2^{(v-1)}\;\beta_{c,v}}  \le \sum\limits_{j\in\mathcal{I}:\exists i, (i,j) \in \mathcal{L}_c}{j \; \psi_{i,j,c}} &\quad& c\in \mathcal{C} \label{logvarmodel_cons2b}\\ 
 &\sum_{(i,j) \in \mathcal{L}_c} \psi_{i,j,c} = 1 &\quad& c\in \mathcal{C}\label{logvarmodel_cons2c}\\
 &  \alpha_{c}, \beta_{c}\in \{0,1\}^{\Theta} &\quad&{c \in \mathcal{C}}\label{logvarmodel_cons3}\\
 & \pi_{u_1,u_2,c,k} \in \{0,1\} &\quad&{c,k\in \mathcal{C}, \, (u_1,u_2) \in \mathcal{N}_{c,k}} \label{logvarmodel_cons4}\\
 & \psi_{i,j,c} \in \{0,1\} &\quad& c \in \mathcal{C} \; (i,j) \in \mathcal{L}_c .  \label{logvarmodel_cons5}
  \end{alignat}
\end{subequations}

In this model, Constraints \eqref{logvarmodel_cons1a}-\eqref{logvarmodel_cons1c} satisfy the non-overlapping feasibility rule, and Constraints \eqref{logvarmodel_cons2a}-\eqref{logvarmodel_cons2c} ensure that each circle is fully contained by the surrounding circle. Finally, Constraints \eqref{logvarmodel_cons3}-\eqref{logvarmodel_cons5} are the domain restrictions for the decision variables. 

\subsection{Relaxed Version}\label{relax_logver}

 In Section \ref{rest_logver}, if a circle is assigned to a candidate point, it means that the center of corresponding circle is located exactly at the corresponding candidate point. In this section, assigning a circle to a candidate point is interpreted differently. In particular, assigning a circle to a candidate point means locating its center to a point included by the region represented with this point. To be more precise, the region represented with a candidate point is the cell whose lower-left corner point is the corresponding candidate point.
The  relaxed version will cover all feasible solutions of the original problem as well as a subset of infeasible solutions; however, the subset of allowed infeasible solutions is limited. 
For example,  locating two circles to two candidate regions is forbidden if even the farthest distance between these two regions is less than the sum of the radii of the corresponding circles.

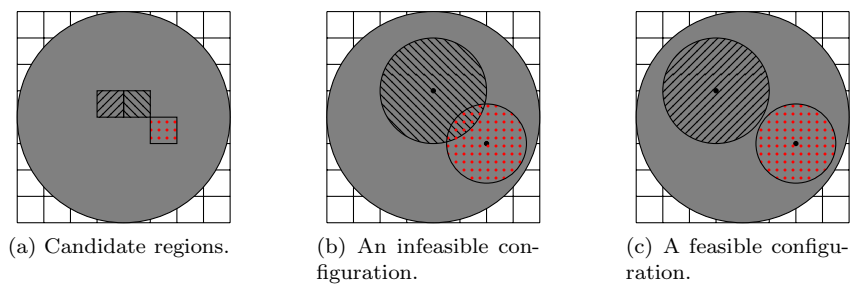
\begin{figure}[H]
\centering
\subfloat[Candidate regions.]{\label{cell_relaxation_example_picture1}
\centering
\begin{tikzpicture}[scale=0.7]
\draw[step=0.5,black,thin, draw opacity = 0.3] (0,0) grid (4,4);

\draw[fill=gray, opacity = 0.15] (2,2) circle (2cm);

\draw[pattern color=black, pattern=north east lines, draw opacity =1, fill opacity = 1] (1.5,2)--(1.5,2.5)--(2,2.5)--(2,2)--cycle;

\draw[pattern color=black, pattern=north west lines, draw opacity = 1, fill opacity = 1] (2,2)--(2,2.5)--(2.5,2.5)--(2.5,2)--cycle;

\draw[pattern color=red, pattern=dots, draw opacity = 1, fill opacity = 1] (2.5,1.5)--(2.5,2)--(3,2)--(3,1.5)--cycle;
\end{tikzpicture}}\hspace{1cm}
\subfloat[An infeasible configuration.]{\label{cell_relaxation_example_picture2}
\centering
\begin{tikzpicture}[scale=0.7]
\draw[step=0.5,black,thin, draw opacity = 0.3] (0,0) grid (4,4);

\draw[fill=gray, opacity = 0.15] (2,2) circle (2cm);

\filldraw[black] (2,2.5) circle (0.04 cm);
\filldraw[black] (3,1.5) circle (0.04 cm);

\filldraw[pattern color=black, pattern=north west lines, draw opacity = 1, fill opacity = 1] (2,2.5) circle (1cm);
\filldraw[pattern color=red, pattern=dots, draw opacity = 1, fill opacity = 1] (3,1.5) circle (0.75cm);
\end{tikzpicture}}\hspace{1cm}
\subfloat[A feasible configuration.]{\label{cell_relaxation_example_picture3}
\centering
\begin{tikzpicture}[scale=0.7]
\draw[step=0.5,black,thin, draw opacity = 0.2] (0,0) grid (4,4);

\draw[fill=gray, opacity = 0.15] (2,2) circle (2cm);

\filldraw[black] (1.5,2.5) circle (0.04 cm);
\filldraw[black] (3,1.5) circle (0.04 cm);

\filldraw[pattern color=black, pattern=north east lines, draw opacity = 1, fill opacity = 1] (1.5,2.5) circle (1cm);
\filldraw[pattern color=red, pattern=dots, draw opacity = 1, fill opacity = 1] (3,1.5) circle (0.75cm);
\end{tikzpicture}}
\caption{Example placements of circles for the relaxed version.}
\label{RelaxedDiscreatization_picture}
\end{figure}

As an example,  look at the candidate regions given in Figure \ref{cell_relaxation_example_picture1}. If  locating the largest circle to the candidate region at the middle and the second largest circle to the candidate region at the rightmost is considered, then    an infeasible configuration is obtained for the proposed relaxed version since the farthest distance of these regions (1.58 units) is smaller than the sum of their radii (1.75 units). However, if the candidate region containing the largest circle's center is changed to the candidate region at the leftmost, then it will be feasible for the relaxed version since the farthest distance between the corresponding regions (1.8 units) is greater than the sum of their radii. This example shows that the   relaxed version will be able to  prevent  some (but not all)  infeasible configurations.

To introduce the relaxed  formulation,  two sets of binary decision variables, $\gamma_c\in\{0,1\}^\Theta$ and $\omega_c\in\{0,1\}^\Theta$, are defined which respectively denote the $x$-axis  and $y$-axis  coordinate of the left-lower corner point of the cell in base 2 where the center of circle $c$ is contained. 
Another variable $\eta_{i,j,c}$ is also defined to ensure that circle $c$ is fully contained by the surrounding circle if there is at least one point of the cell whose left-lower point is contained in the set $ \mathcal{S}_c :=  \left\{i,j \in \mathcal{I}\, : \,  ( i-\theta)^2 + ( j-\theta)^2 \le (\frac{R-r_c}{\delta})^2 \text{ or } (i-\theta+1 )^2 + ( j-\theta+1)^2 \le ( \frac{R-r_c}{\delta})^2 \right\} $.
 
By using a similar idea with the restricted version, another set $\mathcal{O}_{c,k}$ is defined for the pairs of circles $c,k$ similar to $\mathcal{N}_{c,k}$. The set $\mathcal{O}_{c,k}$ includes the pairs of the minimum required number of cells where the farthest distance of the cells between the centers of circles $c,k$ are bigger than the sum of their radii. This set is defined as follows:
\begin{center}
\vspace{0.1cm}
{\small\parbox{14cm}{$\mathcal{O}_{c,k}=\Big\{(u_1,u_2): u_1, u_2 \in\mathcal{I}\backslash\{\theta\},\, u_1,u_2 \ge 0,\,(u_1+1)^2+(u_2+1)^2\ge ((r_c+r_k)/\delta)^2$\par \parindent2.8cm $(u_1\delta)^2+((u_2+1)\delta)^2 < (r_c+r_k)^2,\,((u_1+1)\delta))^2+(u_2\delta)^2 < (r_c+r_k)^2\Big\}$.}}
\vspace{0.1cm}
\end{center}

The defined set $\mathcal{O}_{c,k}$ consists of the minimum number of required cells to locate the circles without overlapping at least if circles $c$ and $k$ are located to the farthest points of corresponding two cells. In other words, the set $\mathcal{O}_{c,k}$ includes the pairs $(u_1,u_2)$ such that the diagonal of the rectangle with side lengths $(u_1+1)$ and $(u_2+1)$ is larger than $(r_c+r_k)/\delta$. Similarly, there are no other rectangles with diagonal length larger than $(r_c+r_k)/\delta$ which is totally included by this rectangle. Also, if centers of circles $c,\,k\in \mathcal{C}$ are included by two cells whose farthest distance is at least $u_1$ cells on $x$-axis and $u_2$ cells on $y$-axis, the corresponding proposed variable $\Pi_{u_1,u_2,c,k}$ will be one. With the help of the sets $\mathcal{O}_{c,k}$ and the variables $\Pi_{u_1,u_2,c,k}$,  the obvious infeasible solutions of CPP are eliminated. 
The corresponding formulation is as follows:
\begin{subequations}\label{logrelaxform}
  \begin{alignat}{3}
 &\sum\limits_{t=1}^{\Theta}{2^{(t-1)}(\gamma_{c,t}-\gamma_{k,t})} \ge  \sum_{(u_1,u_2)\in \mathcal{O}_{c,k}} u_1\; \Pi_{u_1,u_2,c,k} &\quad& c,k\in \mathcal{C}: c<k\label{relaxlogvarmodel_cons1a}\\ 
 &\sum\limits_{v=1}^{\Theta}{2^{(v-1)}(\omega_{c,v}-\omega_{k,v})} \ge  \sum_{(u_1,u_2)\in \mathcal{O}_{c,k}} u_2\; \Pi_{u_1,u_2,c,k} &\quad& c,k\in \mathcal{C}: c<k\label{relaxlogvarmodel_cons1b}\\ 
 &\sum_{(u_1,u_2)\in \mathcal{O}_{c,k}} \Pi_{u_1,u_2,c,k} = 1 &\quad& c,k\in \mathcal{C}: c<k\label{relaxlogvarmodel_cons1c}\\ 
 & \sum\limits_{t=1}^{\Theta}{2^{(t-1)}\;\gamma_{c,t}} \le \sum\limits_{(i,j)\in\mathcal{S}_c}{i \; \eta_{i,j,c}} &\quad& c\in \mathcal{C} \label{relaxlogvarmodel_cons2a}\\  
 &\sum\limits_{v=1}^{\Theta}{2^{(v-1)}\;\omega_{c,v}} \le \sum\limits_{(i,j)\in\mathcal{S}_c}{j \; \eta_{i,j,c}} &\quad& c\in \mathcal{C} \label{relaxlogvarmodel_cons2b}\\ 
 &\sum_{(i,j) \in \mathcal{S}_c} \eta_{i,j,c} = 1 &\quad& c\in \mathcal{C} \label{relaxlogvarmodel_cons2c}\\
 &  \gamma_{c}, \omega_{c}\in \{0,1\}^{\Theta} &\quad&{c \in \mathcal{C}}\label{relaxlogvarmodel_cons3}\\
 & \Pi_{u_1,u_2,c,k} \in \{0,1\} &\quad&{ c,k\in \mathcal{C}, \, (u_1,u_2) \in \mathcal{O}_{c,k}} \label{relaxlogvarmodel_cons4}\\
 & \eta_{i,j,c} \in \{0,1\} &\quad&{\,c \in \mathcal{C}}, \;{ (i,j) \in \mathcal{S}_c}. \label{relaxlogvarmodel_cons5}
  \end{alignat}
\end{subequations}

In this model, Constraints \eqref{relaxlogvarmodel_cons1a}-\eqref{relaxlogvarmodel_cons1c} satisfy the non-overlapping feasibility condition partially while allowing some intersections. Constraints \eqref{relaxlogvarmodel_cons2a}-\eqref{relaxlogvarmodel_cons2c} ensure that there is at least one point of the corresponding cell for locating the circle's center such that it is fully contained by the surrounding circle. Finally, \eqref{relaxlogvarmodel_cons3}-\eqref{relaxlogvarmodel_cons5} are domain restrictions for the decision variables.

\section{Solution Methods and Enhancements}\label{sec:method}

In this section, the solution methods designed for   CPP are presented. A bisection-type solution algorithm is introduced which iteratively considers different candidate radius values for the surrounding circle. The present solution algorithm solves the restricted and relaxed formulations alternately by using the discretized circle while utilizing the ILP models in Section~\ref{sec:ILPformulations}. Then, several algorithmic enhancements are presented, which significantly improve the success of the proposed solution approach.

\subsection{Discretized-Space Circle Packer Algorithm}\label{sec_solutionmethods}

In this subsection,  the details of the solution procedure, named as the Discretized-Space Circle Packer Algorithm (DCPACK), are given. This algorithm depends on updating the upper and lower bounds of the radius of the surrounding circle until they are close enough.  
 To initialize these bounds, it is  proceeded as follows: If  a surrounding circle in which all circles can be packed is known, then its radius gives an upper bound. Otherwise, the upper bound is trivially selected as the sum of the radii of all circles. Similarly, if  a surrounding circle into which there is no feasible placement of packing all circles is known, then its radius gives a lower bound for the optimal solution; otherwise, the lower bound is initialized as the sum of the radii of the largest two circles. 

\begin{algorithm}
\caption{DCPACK Algorithm.}	\label{alg:circlePacking}
		\begin{algorithmic}[1]
		    \REQUIRE $\mathcal{C}$; $r_c, \; \forall c\in \mathcal{C}$; $U$; $L$; $\delta$; $\varepsilon$.  
            \ENSURE $(x_c,y_c),\; \forall c\in \mathcal{C}$; $U$; $L$.
            \STATE $R \leftarrow \frac{U+L}{2}$ \label{CPA_Step1}
			\STATE Divide the container with radius $R$ into square cells with length   $\delta$. \label{CPA_Step2}
			\STATE Solve the restriction model~\eqref{logrestform} \revise{using an MILP solver}. \label{CPA_Step3}
            \IF {the restriction model is feasible}
			\STATE $U \leftarrow R$. Go to Step~\ref{CPA_StepRecurrence}.\label{CPA_Step6}
			\ENDIF
			\STATE Solve the relaxation model~\eqref{logrelaxform} \revise{using an MILP solver}. \label{CPA_Step9} 
            \IF {the relaxation model is infeasible}
			\STATE $L \leftarrow R$. Go to Step~\ref{CPA_StepRecurrence}.\label{CPA_Step13}
			\ELSE
			\STATE $\delta \leftarrow \frac{\delta}{2}$. Go to Step~\ref{CPA_Step2}.\label{CPA_Step11}
			\ENDIF
			\IF {$U-L > \varepsilon U$} \label{CPA_StepRecurrence}
			\STATE Go to Step~\ref{CPA_Step1}.
			\ELSE
			\STATE STOP!
			\ENDIF
		\end{algorithmic}
\end{algorithm}

 Algorithm~\ref{alg:circlePacking}  aims to progressively improve the upper and lower bounds for the radius of the surrounding circle \revise{and terminates when the $\varepsilon$-optimality is proven, which is defined as follows: Suppose  there is   a feasible configuration for which the radius of the container is $U$ and it is proven that there does not exist any feasible configuration when the radius of the container is at most $L$. Then,  $U$ is called   an $\varepsilon$-optimal solution for the radius of the container if condition $U-L > \varepsilon U$ holds.}

First, the surrounding circle's radius is initialized as the average of the initial upper and lower bounds at Step~\ref{CPA_Step1}. Then, the given surrounding circle is divided into smaller squares with side length $\delta$ at Step~\ref{CPA_Step2}  and   the restricted  formulation~\eqref{logrestform} is solved. If this problem gives a feasible configuration, the corresponding radius value gives an upper bound for the minimum value of the surrounding circle's radius. Then,   the upper bound is updated at Step~\ref{CPA_Step6}, and  the gap between the bounds are is checked at Step~\ref{CPA_StepRecurrence}. If the bounds are close enough, then Algorithm~\ref{alg:circlePacking} terminates with an $\varepsilon$-optimal solution. Otherwise, Algorithm~\ref{alg:circlePacking} proceeds with Step~\ref{CPA_Step1} with an updated upper  bound.

On the other hand, if the restricted model is infeasible,   the relaxed  formulation~\eqref{logrelaxform} is solved.   If the relaxed model is also infeasible for the given radius of the surrounding circle, it means that there is no feasible configuration of the given circles into this surrounding circle. Hence,    the lower bound for the original problem is updated at Step~\ref{CPA_Step13} and it is  checked whether the updated lower and upper bounds are close enough or not. If they are close enough, the algorithm stops. Otherwise, Algorithm~\ref{alg:circlePacking} proceeds with Step~\ref{CPA_Step1} with an updated lower  bound.

If a given radius value for the surrounding circle is feasible in the relaxation model which is infeasible in the restriction model, then there is no clear conclusion. In that case,   the side length of the cells is halved to get a finer discretization at Step~\ref{CPA_Step11} and  the algorithm proceeds with  the new discretization.  

\revise{Now, the convergence of Algorithm~\ref{alg:circlePacking} is discussed. The algorithm will set the value of $R$ at most $\lceil \log_2(1/\varepsilon) \rceil$ many times in Step~\ref{CPA_Step1}, a property inherited from the bisection method. As long as the algorithm executes Step~\ref{CPA_Step11} finitely many times,   the finite convergence to an $\varepsilon$-optimal solution can be shown. The only complication may arise when $L$ and $U$ are not close enough, but the current iterate $R$ is arbitrarily close to the global optimal value. In this case, the restriction model~\eqref{logrestform} will be infeasible and the relaxation model~\eqref{logrelaxform} will be feasible, and the algorithm might have to refine the grid arbitrarily many times in Step~\ref{CPA_Step11}. Although  this behavior has not observed in the computational experiments, this occurrence can be prevented by additional measures. For example, if the algorithm is stuck at Step~\ref{CPA_Step11} for a predetermined number of iterations, then the current iterate $R$   can be perturbed  slightly and the algorithm can continue with Step~\ref{CPA_Step2}.}

\subsection{Algorithmic Enhancements}\label{sec_algorithmicenhancements}

In this subsection,   the details of  improvements proposed for Algorithm~\ref{alg:circlePacking} are given.  In Section \ref{sec:solutionspacereductions},  the details of Algorithm~\ref{alg:identifyingFeasibleRegions}, which decreases the number of decision variables by decreasing the solution space for each circle, are presented. In Section \ref{sec:lowerbound},  the methods for obtaining initial lower \revise{and upper} bounds for the radius of the surrounding circle are described used in Algorithm~\ref{alg:circlePacking}. Then, other improvements are given in Section \ref{sec:otherimp}.
Throughout this subsection,  the circles are assumed to be ordered with respect to their radii, that is, $r_c \ge r_k$ for $c <k$.

\subsubsection{Solution Space Reductions}\label{sec:solutionspacereductions}  In this part, the aim is to reduce the feasible region of CPP using geometric arguments. 
 For example, without loss of generality, it is possible to place the center of the first circle to the first quadrant, and the center of the second circle  to the half space defined by $\{(x,y)\in \mathbb{R}^2\, : \; y-x \ge 0\}$ by adding the constraints $x_1 \ge 0, y_1 \ge 0, y_2 \ge x_2$. This helps  eliminating some symmetric solutions of the problem.

In fact, more restrictions   can be inferred  by considering the packing and non-overlapping   constraints. To start with,  it is possible to use the fact that that the distance from a point to the boundary of the surrounding circle should be at least the radius of the circle to be located for placing at this point; and the maximum distance from the point to the boundary should be at least two times of the radius of the reference circle (the largest circle for other circles, and the second largest circle for the largest circle) plus the radius of the corresponding circle by considering the non-overlapping constraint. Moreover,   the obtained regions    can  be further reduced  by eliminating the points where there is no feasible point to locate the other circle's center due to the non-overlapping constraint.  These resulting regions are the smallest feasible regions for these circles within this procedure. Then, by using the reduced regions  for the largest two circles, the feasible regions for locating centers of other circles are reduced  by performing an iterative algorithm.  The steps of this approach is given in Algorithm~\ref{alg:identifyingFeasibleRegions}

\revise{
An example is given in Figure \ref{solspacered_picture} to illustrate the solution space reduction idea. Suppose that in a given instance, two circles have radii 6 and 7 units, and the task it to identify the regions their centers can be located to a container of radius 13.6 units. 
In Figure \ref{redregver1}, the dashed and   undashed regions show the reduced regions to place the largest and the second largest circles, respectively. 
Then, the improved version of the reduced regions are given in Figure \ref{redregver1imp1}. 
These reduced regions given in Figure \ref{redregver1imp1} will be combined with the knowledge that any two circles cannot overlap. Then, the corresponding regions can be decreased by eliminating the points where there is no feasible point to locate the other circle's center. After this elimination step,   feasible regions given in Figure \ref{redregver1imp2} are obtained for the largest two circles. These resulting regions are the minimum possible feasible regions for these circles within this procedure.

\begin{figure}[H]
\centering
\subfloat[Initial reduced regions.]{\label{redregver1}
\centering
\begin{tikzpicture}[scale=0.20]
\draw[draw opacity=0.0] (-8.8,0) rectangle (8.4,0);

\filldraw[gray, fill opacity=0.1] (0,0) circle (6.8 cm);
\draw[red, draw opacity=0.2] (-4.8,-4.8) -- (4.8,4.8);

\draw[fill=red, opacity=0.05] (0,0) -- +(45:6.8cm)  arc (45:225:6.8cm) -- cycle;
\draw[fill=blue, opacity=0.05] (0,0) -- +(0:6.8cm)  arc (0:90:6.8cm) -- cycle;
\draw[pattern=north west lines, pattern color=blue, opacity=0.4] (0,0) -- +(0:6.8cm)  arc (0:90:6.8cm) -- cycle;

\end{tikzpicture}
}\hspace{0.1cm}
\subfloat[Improved reduced regions.]{\label{redregver1imp1}
\centering
\begin{tikzpicture}[scale=0.20]
\draw[draw opacity=0.0] (-9.2,0) rectangle (9,0);

\filldraw[gray, fill opacity=0.1] (0,0) circle (6.8 cm);

\fill[red, opacity=0.1, even odd rule] (0,0) -- +(45:3.8cm)  arc (45:225:3.8cm) -- cycle (0,0) -- +(45:3.2cm)  arc (45:225:3.2cm) -- cycle;

\fill[blue, opacity=0.1, even odd rule] (0,0) -- +(0:3.3cm)  arc (0:90:3.3cm) -- cycle (0,0) -- +(0:2.7cm)  arc (0:90:2.7cm) -- cycle;

\fill[pattern=north west lines, pattern color=blue, opacity=0.6, even odd rule] (0,0) -- +(0:3.3cm)  arc (0:90:3.3cm) -- cycle (0,0) -- +(0:2.7cm)  arc (0:90:2.7cm) -- cycle;

\end{tikzpicture}
}\hspace{0.1cm}
\subfloat[The resulting reduced regions.]{\label{redregver1imp2}\centering\begin{tikzpicture}[scale=0.20]
\draw[draw opacity=0.0] (-10.3,0) rectangle (10.3,0);

\filldraw[gray, fill opacity=0.1] (0,0) circle (6.8 cm);

\coordinate (A) at (225:6.8);
\draw[blue] (0,0) +(0:2.7cm) arc (0:45:2.7cm);
\draw[blue] (A) +(45:9.5cm) arc (45:58:9.5cm);
\draw[blue] (3.3,0) arc (0:87:3.3cm);
\draw[blue] (2.7,0) -- (3.3,0);

\draw[red, rotate=180] (3.2,0) arc (0:45:3.2cm);
\draw[red, rotate=225] (3.2,0) -- (3.8,0);
\draw[red, rotate=132] (3.8,0) arc (0:93:3.8cm);
\draw[red] (-3.2,0) arc (180:154:6.5cm);
\end{tikzpicture}}
\caption{An example of solution space reductions.}
\label{solspacered_picture}
\end{figure}
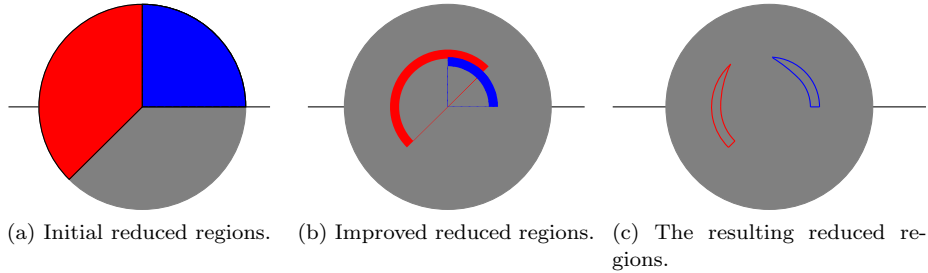

}

\begin{algorithm}
\caption{Feasible Region Identification Algorithm.}	\label{alg:identifyingFeasibleRegions}
		\begin{algorithmic}[1]
		    \REQUIRE $\mathcal{C}$; $r_c, \; \forall c\in \mathcal{C}$; $\mathcal{S}_1$; $\mathcal{S}_2$; $R$, $k$.  
            \ENSURE True or False.
            \FOR {$\{p_c=(x_c,y_c): \lVert (p_c,p_l)\rVert_2 \ge (r_c+r_l), \lVert p_c\rVert_2 \le R-r_c\}$ where the center of circle at order $c$ is located where $c<k$}
                \STATE $\mathcal{G} = \{p_k=(x_k,y_k)\in \mathbb{R}^2: (R-r_k) \ge \lVert p_k\rVert_2\}\ge \max\{0, 2 \, r_1 + r_k - R\}$
                \STATE $\mathcal{F} = \left\{p_k=(x_k,y_k)\in\mathbb{R}^2: \forall c\in\mathcal{C}, \lVert (p_k,p_c) \rVert_2\ge r_c+r_k\right\}\cap\mathcal{G}$\label{IFRA_Step0}
            \ENDFOR
            \IF {$\mathcal{F} = \emptyset$}
                \STATE return FALSE.
            \ELSE
                \STATE return TRUE.
            \ENDIF
		\end{algorithmic}
\end{algorithm}

It is easy to see that if Algorithm~\ref{alg:identifyingFeasibleRegions} results in False for any $k$, there is no feasible configuration of circles in the given surrounding circle; otherwise, there is at least one feasible placement. By using the resulting reduced regions obtained by Algorithm~\ref{alg:identifyingFeasibleRegions}, it is possible to eliminate the decision variables corresponding to the candidate points which are not included by the resulting regions for each circle.

\subsubsection{Initializing Lower \revise{and Upper} Bounds}\label{sec:lowerbound}

\revise{An obvious but potentially weak upper bound can be obtained as the sum of all radii of the circles, that is, $\sum_{c\in\mathcal{C}} r_c$. 
To obtain better upper bound values,} it is possible to use the best-known upper bounds from Circle Packing Contest of Al Zimmermann's Programming Contests and Packomania website with the algorithm given in \cite{minRadius_Huang2006}. 

Below, four methods to  initialize the lower bound (LB) in Algorithm~\ref{alg:circlePacking}  will be discussed. Many of these methods exploit the geometric properties of the problem whereas the last one solves an auxiliary optimization problem to compute the bound. 

\begin{itemize}
\item LB1: An obvious lower bound for the radius of the surrounding circle is the sum of  the radii of the largest two circles, that is, $r_1+r_2$.  
 \item LB2: Since the sum of the areas of all circles gives  a lower bound for the area of the surrounding circle,   another lower bound for  the radius of the surrounding circle can be obtained  as  $\sqrt{\sum_{c\in\mathcal{C}}r_c^2}$.
\item LB3:  Algorithm~\ref{alg:identifyingFeasibleRegions} can be performed iteratively to calculate another lower bound for  the radius of the surrounding circle. 
\item LB4: Placing circles on the plane, even in tangent positions, inevitably creates idle regions that cannot contain any other circle. The main idea in this lower bound calculation is to quantify the area of idle regions between circle triples and  and the circles adjacent to the surrounding circle. This approach is summarized below and  the reader is referred to  \cite{tacspinar2021discretization} for details.

\revise{
To start with, define a set denoted  by $\overline{\mathcal{C}}=\mathcal{C}\cup\{0\}$, where 0 is the index of the surrounding circle.
In order to   investigate the idle regions, it is convenient to think of  a configuration of circles as a graph such that each circle denotes a node, and there is an edge between any two circles if their centers can be connected by a line without intersecting   another circle. As an illustration,  consider Figure~\ref{adjacentcirclesfigure}. Observe that there are three adjacent circles in Figure~\ref{lowerbound_threeadjcircles} and they are all tangent to one another. In general, it is possible to compute  the idle area between such triplet of circles $C$, $K$ and $L$, denoted by $\Delta_{C,K,L}$, as follows:
\[
{\sqrt{(r_C+r_K+r_L) r_C r_K r_L } \,- \, \sum\limits_{\substack{ (c,k,l)\in\{(C,K,L),\\(K,C,L),(L,C,K)\}}}{\frac{r_c^2 \cos^{-1}\left(\frac{r_c(r_c+r_k+r_l)-r_kr_l}{(r_c+r_l)(r_c+r_k)}\right)}{2}}}.
\]
It can also happen that the circles are not tangent to each other, or form cycles with length more than three as in Figure~\ref{4adjacentcirclesfigure1}. In these cases, it is possible to underestimate the idle area  as the sum of the idles areas of the consecutive triples as if they are    adjacent to one another. 
}

\begin{figure}[H]
\centering
\subfloat[Three adjacent circles.]{\label{lowerbound_threeadjcircles}\centering\begin{tikzpicture}[scale=0.25]
\filldraw[draw opacity=0.0, fill opacity=0.0] (-6.2,-1) rectangle (9.5,1);
\filldraw[gray, fill opacity=0.55] (1.3,0.65) circle (1.8 cm);
\filldraw[gray, fill opacity=0.30] (-0.5,4) circle (2 cm);
\filldraw[gray, fill opacity=0.05] (4,4) circle (2.5 cm);
\node[text width=0.2cm] at (1.3,0.65) { };
\node[text width=0.2cm] at (-0.5,4) { };
\node[text width=0.2cm] at (4,4) { };
\end{tikzpicture}
}\hspace{0.5cm}
\subfloat[Four adjacent circles.]{\label{4adjacentcirclesfigure1}\centering\begin{tikzpicture}[scale=0.25]
\filldraw[draw opacity=0.0, fill opacity=0.0] (-9,-1) rectangle (12,1);
\node[text width=0.2cm] at (1.3,0.65) { };
\node[text width=0.2cm] at (-3.4,4) { };
\node[text width=0.2cm] at (3.2,6.8) { };
\node[text width=0.2cm] at (6.55,3.15) { };
\filldraw[gray, fill opacity=0.55] (1.3,0.65) circle (1.8 cm);
\filldraw[gray, fill opacity=0.30] (-3.4,4) circle (2 cm);
\filldraw[gray, fill opacity=0.05] (3.2,6.8) circle (2.5 cm);
\filldraw[gray, fill opacity=0.15] (6.55,3.15) circle (2.2 cm);
\draw[draw opacity=0.2] (1.3,0.65) -- (-3.8,4) -- (3.2,6.8) -- (6.55,3.15) -- cycle;
\end{tikzpicture}}\hspace{0.5cm}
\caption{The idle region between adjacent circles.}
\label{adjacentcirclesfigure}
\end{figure}
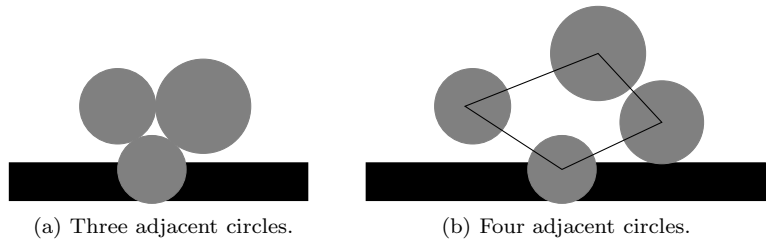

\revise{
There are  other possible configurations  as illustrated in Figure~\ref{lowerbound_surroundingadjacentcircles}. As in Figure~\ref{lowerbound_surroundingadjacent}, if one of the circles in the triplet is the surrounding circle, then the idle area calculation  changes  and  involves an upper bound on the radius of the container. Figure~\ref{lowerbound_surroundingadjacent_3} presents a case in which some idle regions are double counted.  If   circle $k$ is the circle between the circles $c$, $l$ and the surrounding circle, then  the double counted area is computed as 
$\rho_{c,k,l} :=\Delta_{0,c,k}+\Delta_{0,k,l}+\Delta_{c,k,l}+\pi r_k^2$.
}

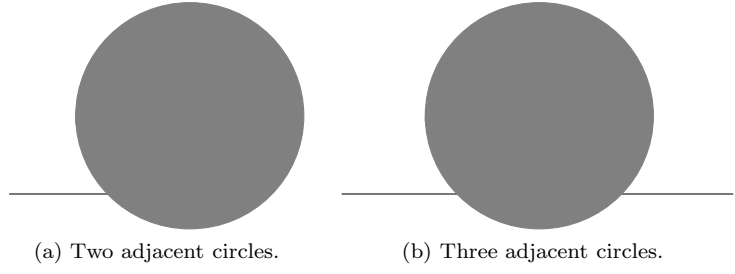
\begin{figure}[H]
\centering
\subfloat[Two adjacent circles.]{\label{lowerbound_surroundingadjacent}
\centering
\begin{tikzpicture}[scale=0.25]
\draw[draw opacity=0.0] (-6.0,0) rectangle (6.9,0);

\filldraw[gray,  fill opacity=0.55] (1.3,0.65) circle (1.8 cm);
\filldraw[gray,  fill opacity=0.30] (-0.5,4) circle (2 cm);
\filldraw[gray, fill opacity=0.05] (3.5,4.15) circle (6.0 cm);
\end{tikzpicture}
}\hspace{0.1cm}
\subfloat[Three adjacent circles.]{\label{lowerbound_surroundingadjacent_3}
\centering
\begin{tikzpicture}[scale=0.25]
\draw[draw opacity=0.0] (-6.9,0) rectangle (13.7,0);

\filldraw[gray,  fill opacity=0.55] (1.3,0.65) circle (1.8 cm);
\filldraw[gray, fill opacity=0.30] (-0.5,4) circle (2 cm);
\filldraw[gray, fill opacity=0.95] (-0.85,1.47) circle (0.5 cm);
\filldraw[gray, fill opacity=0.05] (3.5,4.15) circle (6.0 cm);
\end{tikzpicture}}
\caption{The idle region between adjacent circles and the surrounding circle.}
\label{lowerbound_surroundingadjacentcircles}
\end{figure}

\revise{
It is now explained how to use an integer programming formulation
to calculate a lower bound for the idle area in the optimal packing. In addition to the idle area parameters  $\Delta_{c,k,l}$ and $\rho_{c,k,l}$ as explained above,    also define parameters  $\underline\kappa_{c}$ and $\overline\kappa_c$, which respectively give  the minimum and maximum number of triples circle $c$ can have as neighbors. These parameters can be easily computed as a function of the radius $r_c$.
}
Two binary variables $d_{c,k,l}$ and $f_{c,k,l}$ are defined for  $c,k,l\in\overline{\mathcal{C}}$ with $ c < k < l$. If the circles $c$, $k,$ and $l$ are adjacent to each other, then  $d_{c,k,l}$ takes value 1.  If circles $c$, $k$ and $l$ are adjacent to the surrounding circle    (i.e., $d_{c,k,l}=1$, $d_{0,c,k}=1$, $d_{0,k,l}=1$, $d_{0,c,l}=1$), then the binary variable $f_{c,k,l}$ takes value 1.  The resulting integer linear program is given as below:
\begin{subequations}\label{integerprogram}
  \begin{alignat}{3}
 \qquad\min &\; \sum\limits_{c,k,l\in\overline{\mathcal{C}}: c<k<l}\left({\Delta_{c,k,l}d_{c,k,l}}-{\rho_{c,k,l}f_{c,k,l}}\right) &&\label{integerprogram_obj}\\ 
 \text{s.t.}  
 \  & \underline\kappa \le  \sum\limits_{k,l\in\overline{\mathcal{C}}}d_{c,k,l} + \sum\limits_{k,l\in\overline{\mathcal{C}}}d_{k,c,l} + \sum\limits_{k,l\in\overline{\mathcal{C}}}d_{k,l,c} \le \overline\kappa &\qquad&{c \in \overline{\mathcal{C}} } \label{integerprogram_cons2}\\
 \  & f_{c,k,l} \le d_{c,k,l} &\qquad&{ c,k,l \in \overline{\mathcal{C}} } \label{integerprogram_cons3}\\
 \  & f_{c,k,l} \le d_{0,c,k},\quad f_{c,k,l} \le d_{0,k,l},\quad  f_{c,k,l} \le d_{0,c,l}  &\qquad&{ c,k,l \in \overline{\mathcal{C}} } \label{integerprogram_cons4}\\
 \  & f_{c,k,l} \ge d_{c,k,l} + d_{0,c,k} + d_{0,k,l} + d_{0,c,l} - 3 &\qquad&{ c,k,l \in \overline{\mathcal{C}} } \label{integerprogram_cons5}\\
 \  & d_{c,k,l}, f_{c,k,l} \in \{0,1\} &\qquad&c,k,l \in \overline{\mathcal{C}} . \label{integerprogram_cons6}
  \end{alignat}
\end{subequations}
\revise{
In the above formulation, objective function~\eqref{integerprogram_obj} minimizes a lower bound on the total idle area created by circle triplets. Constraint~\eqref{integerprogram_cons2} ensures that circle $c$ has at least 
$\underline\kappa_{c}$ and at most $\underline\kappa_c$ many
circles as its neighbors.  
 Constraints~\eqref{integerprogram_cons3}-\eqref{integerprogram_cons5}   assign the value of the decision
variable $f_{c,k,l}$ according to the values of the decision variables $d_{c,k,l}$, $d_{0,c,k}$, $d_{0,k,l}$ and $d_{0,c,l}$. In particular, $f_{c,k,l}$ takes value one if and only of    $d_{c,k,l}$, $d_{0,c,k}$, $d_{0,k,l}$ and $d_{0,c,l}$ are all assigned to one. Finally, Constraint~\eqref{integerprogram_cons6} are  domain restrictions for the decision variables.
}
The sum of the optimal value of problem~\eqref{integerprogram} and the sum of the areas of all circles gives    a lower bound for the area of the surrounding circle, from which  a   lower bound for  the radius of the surrounding circle   can be computed.
\end{itemize}

\subsubsection{Other Improvements}\label{sec:otherimp}

In addition to the previous enhancements,  a set of procedures handled during
the branching process are also introduced.  First of all, the decision variables are prioritized according
to the radius of the corresponding circle. The variables corresponding larger circles
have the larger priorities. \revise{This is implemented by Gurobi solver's branching priority functionality}.

\revise{In a typical node of a branch-and-bound tree, some binary variables are assigned a value of one. In the case of CPP, this means that a  number of   circles are assigned to certain candidate points or cells, and some circles are currently unassigned. 
During the branching process, it is possible to encounter situations in which a subset of unassigned circles cannot be located due to insufficient space left. For example, this can happen   if there is not enough area to locate some of the unassigned circles, or  there is not  enough candidate points
to place all the unassigned circles. 
As listed below,   these cases are identified during the execution of an MILP solve, and   cuts are added or    branch-and-bound nodes  are eliminated through   Gurobi solver's user cut and callback functionalities to handle such situations efficiently:
\begin{itemize}
\item The idle regions between the assigned circles where none of the unassigned
circles can be placed are identified. The candidate points or cells in such regions are eliminated.

\item The remaining idle area is calculated  and it is compared with the sum of the areas of the unassigned circles with the minimum
idle area between them to be packed. A branch-and-bound node is eliminated in two cases: i) if the remaining area is less than the sum
of the areas of the unassigned circles, or ii)  if the distance between the farthest available candidate points  is less
than the sum of the radii of the two largest unassigned circles.

\item Situations in which the assignment of an unassigned circle to a particular cell creates overlapping with the already assigned circles are also identified. Then,   such cases are prevented with the help of conditional constraints.
\end{itemize}
}

\section{Computational Study} \label{sec:computational}

In this section, the present algorithm is evaluated on a series of benchmark instances    from \cite{packomania}, the Circle Packing Contest of Al Zimmermann's Programming Contests from \cite{zimmcontest}, and   equal circle instances from \cite{GlobalSolver_HUANG2011}. The instances from the Circle Packing Contest are denoted with $``Zimm-n"$, which represents an instance with $n$-circles where $5\le n\le 20$ \citep{zimmcontest}, and  $n$   equal circle instances from \cite{GlobalSolver_HUANG2011} are denoted with $``Eq-n"$. 

Before analyzing the performance of the algorithm,  the problem~\eqref{quadForm} is solved with Baron and Gurobi  to show the effectiveness of the global solvers in Subsection \ref{compSt:BaronGurobi}. Then,  the performance of the discretization-based solution approach is given in Subsection \ref{sec:performance of discr}.  
The present algorithm for the equal circle case is compared with a heuristic method proposed in \cite{GlobalSolver_HUANG2011} in Subsection \ref{compSt:EqualCase}. Then, the present solution procedure is compared with another heuristic designed for a  rectangular-shaped container from \cite{GlobalSolver_Stoyan2008} in Subsection \ref{compSt:UnequalCase}.   The optimality gap is set as 1\% for Baron and Gurobi  as well as Algorithm~\ref{alg:circlePacking}.

\subsection{Performance Analysis of Global Solvers}\label{compSt:BaronGurobi}

In this subsection,      four instances from \citep{zimmcontest} are used to analyze the performance of global solvers for two versions of formulation~\eqref{quadForm}. For the basic version, no knowledge is shared with the global solvers in terms of solution space reductions \revise{from Section~\ref{sec:solutionspacereductions} or lower and upper bounds from Section~\ref{sec:lowerbound}}  whereas  stronger lower and upper bounds are used for the improved version. These formulations are implemented in Python, and solved   by solvers Baron and Gurobi   with a time limit of 30 minutes. The results  are reported in 
Table~\ref{Tab:BaronGurobiResults}.
\begin{table}[H]
		\caption{Results obtained by the global solvers (TLE: time limit exceeded).} \label{Tab:BaronGurobiResults}
	\begin{center}
    \resizebox{\textwidth}{!}{
		\begin{tabular}{|c|ccc|ccc|ccc|ccc|}
			\hline
			&\multicolumn{6} {c|} {\bfseries Baron} & \multicolumn{6} {c|} {\bfseries Gurobi }\\
     		& \multicolumn{3} {c|} {\underline{Basic Version}} & \multicolumn{3} {c|} {\underline{Improved Version}} &\multicolumn{3} {c|} {\underline{Basic Version}} & \multicolumn{3} {c|} {\underline{Improved Version}} \\ 
     		Ins. &	\shortstack{Upper\\Bound} &	\shortstack{Solution\\ Time (m)}	&	\shortstack{GAP (\%)}	&	\shortstack{Upper\\Bound} &	\shortstack{Solution\\ Time (m)}	&	\shortstack{GAP (\%)}	&	\shortstack{Upper\\Bound} &	\shortstack{Solution\\ Time (m)}	&	\shortstack{GAP (\%)} & \shortstack{Upper\\Bound} &	\shortstack{Solution\\ Time (m)}	&	\shortstack{GAP (\%)}	\\  \hline
$Zimm-5$ & 9.001 & 1.26 & 0.01 & 9.001 & 0.01 & 0.00 & 9.001 & 0.50 & 0.99 & 9.034 & 0.72 & 0.36\\ 
$Zimm-6$ & 11.057 &  12.83 &	0.01 & 11.109  &	0.27 &	0.46 & 11.057 &  17.31 &	0.52 & 11.062  &	1.28 &	0.04\\
$Zimm-7$ & 13.462 & TLE & 3.43 & 13.462 &	TLE	&	2.28 & 13.463 & TLE & 3.44  &	13.462	&	TLE	&	2.28\\ 
$Zimm-8$ & 16.222 & TLE & 7.53 & 16.222 &	TLE	& 5.85	 & 16.384 & TLE & 8.45 &  16.437 &	TLE	& 7.08\\ 
			\hline
		\end{tabular}
        }
	\end{center}
\end{table} 
It is  observed that the global solvers cannot solve problem~\eqref{quadForm} even for the 7-circle
instance for both versions. In the improved version, the given
lower bounds are not improved during the solution times by the solvers, but the
global solvers can detect and improve upper bounds for problem~\eqref{quadForm}. It is deduced that CPP is quite challenging to solve with the help of the global solvers even for
the small-sized instances.

\subsection{Performance Analysis of the Discretization-Based Approach}
 \label{sec:performance of discr}

In this subsection, the performance of Algorithm~\ref{alg:circlePacking} and the effects of the enhancements given in Section \ref{sec_algorithmicenhancements} are   investigated. In Table~\ref{tab:algenhancements},  the results obtained by Algorithm~\ref{alg:circlePacking} are given for two different versions: the basic version without any  enhancements and the enhanced version with all the enhancements proposed in Section \ref{sec_algorithmicenhancements}. As seen in Table~\ref{tab:algenhancements}, the enhancements improve Algorithm~\ref{alg:circlePacking}'s performance, and they decrease the total solution time significantly.  
\begin{table}[H]
	\vspace{0.5cm}
		\caption{Performance analysis of the enhancements introduced for Algorithm~\ref{alg:circlePacking}.} \label{tab:algenhancements} 
	\begin{center}
    \resizebox{\textwidth}{!}{
		\begin{tabular}{|c|ccc|cccc|}
     		\hline
     		& \multicolumn{3} {c|} {\bfseries Basic Version} & \multicolumn{4} {c|} {\bfseries Enhanced Version} \\ 
     		Instance	&	\shortstack{Upper\\Bound}&	\shortstack{Total Sol.\\Time (m)}	&	\shortstack{GAP (\%)}	&	\shortstack{Upper\\Bound} &	\shortstack{Pre-processing \\Time (m)} &	\shortstack{Total Sol.\\Time (m)}	&	\shortstack{GAP (\%)}	\\  \hline
			$Zimm-5$ & 9.002 & 2.07 & 0.02 & 9.001 & 0.04 & 0.06 & 0.76 \\
			$Zimm-6$ & 11.139 & 9.25 & 0.73 & 11.071 & 0.17 & 0.24 & 0.13 \\
			$Zimm-7$ & 13.471 & 17.72 & 0.82& 13.467 & 0.86 & 1.07 & 0.64 \\
			$Zimm-8$ & 16.321 & 28.20 & 0.99 & 16.224 & 1.59 & 2.57 & 0.81 \\
            \hline
		\end{tabular}
        }
	\end{center}
\end{table} 
The individual effect of each enhancement is also analyzed. For this purpose,   four versions of Algorithm~\ref{alg:circlePacking} are tested \revise{in which all enhancements from 
Section~\ref{sec_algorithmicenhancements} are implemented similar to Enhanced Version, except for one enhancement left out each time. In particular, Version 1 is the one without any lower bound improvement 
 from Section~\ref{sec:lowerbound}, Version 2 is the one without any upper bound improvement 
from Section~\ref{sec:lowerbound}, Version 3 is the one without any solution space reductions 
from Section~\ref{sec:solutionspacereductions}, and Version 4 is the one without any cuts added during the branching process from Section~\ref{sec:otherimp}. The results are presented   in Tables~\ref{tab:enhancementsonebyone1} and  \ref{tab:enhancementsonebyone2}.}
\begin{table}[H]
	\vspace{0.5cm}
		\caption{Effects of the qualities of starting lower and upper bounds for Algorithm~\ref{alg:circlePacking}.} \label{tab:enhancementsonebyone1} 
	\begin{center}
    \resizebox{\textwidth}{!}{
		\begin{tabular}{|c|cccc|cccc|cccc|cccc|}
     		\hline
     		& \multicolumn{4} {c|} {\bfseries Version 1} & \multicolumn{4} {c|} {\bfseries Version 2}  \\ 
     		Instance	&	\shortstack{Upper\\Bound} &	\shortstack{Pre-processing \\Time (m)} &	\shortstack{Total Sol.\\Time (m)}	&	\shortstack{GAP (\%)} &	\shortstack{Upper\\Bound} &	\shortstack{Pre-processing \\Time (m)} &	\shortstack{Total Sol.\\Time (m)}	&	\shortstack{GAP (\%)} \\  \hline
			$Zimm-5$ & 9.002 & 0.01 & 1.43 & 0.02 & 9.006 & 0.03 & 2.12 & 0.07 \\
			$Zimm-6$ & 11.140 & 0.04 & 6.63 & 0.86 & 11.084 & 0.14 & 5.73 & 0.43 \\
			$Zimm-7$ & 13.469 & 0.15 & 13.56 & 0.89 & 13.470 & 0.68 & 11.41 & 0.80 \\
			$Zimm-8$ & 16.273 & 0.41 & 21.32 & 0.88 & 16.269 & 1.11 & 18.83 & 0.83 \\
			\hline
		\end{tabular}
        }
	\end{center}
\end{table}
Initializing lower and upper bounds is an important factor in Algorithm~\ref{alg:circlePacking} according to Table~\ref{tab:enhancementsonebyone1}. The pre-processing time required for initializing the bounds is less than the necessary time spent in Algorithm~\ref{alg:circlePacking} without the pre-processing methods. Hence, the methods for initializing the upper and lower bounds are effective since the solution times do not decrease much compared to the Basic Version.
\begin{table}[H]
	\vspace{0.5cm}
		\caption{Effects of other enhancements one by one introduced for Algorithm~\ref{alg:circlePacking}.} \label{tab:enhancementsonebyone2} 
	\begin{center}
    \resizebox{\textwidth}{!}{
		\begin{tabular}{|c|cccc|cccc|cccc|cccc|}
     		\hline
     		& \multicolumn{4} {c|} {\bfseries Version 3} & \multicolumn{4} {c|} {\bfseries Version 4}  \\ 
     		Instance	&	\shortstack{Upper\\Bound} &	\shortstack{Pre-processing \\Time (m)} &	\shortstack{Total Sol.\\Time (m)}	&	\shortstack{GAP (\%)} &	\shortstack{Upper\\Bound} &	\shortstack{Pre-processing \\Time (m)} &	\shortstack{Total Sol.\\Time (m)}	&	\shortstack{GAP (\%)} \\  \hline
			$Zimm-5$ & 9.001 & 0.03 & 0.06 & 0.86 & 9.001 & 0.04 & 0.07 & 0.55 \\ 
			$Zimm-6$ & 11.091 & 0.15 & 0.26 & 0.36 & 11.087 & 0.17 & 0.28 & 0.25 \\
			$Zimm-7$ & 13.466 & 0.81 & 1.04 & 0.70 & 13.468 & 0.84 & 1.13 & 0.59 \\
			$Zimm-8$ & 16.224 & 1.53 & 6.48 & 0.69 & 16.224 & 1.56 & 2.91 & 0.86 \\
			\hline
		\end{tabular}
        }
	\end{center}
\end{table} 
According to Table~\ref{tab:enhancementsonebyone2}, one can say that the solution space reductions and adding cuts during the branching have a positive effect in the reduction of the total solution time compared to the Basic Version. However, they are less influential compared to the lower and bound enhacement methods.

Other than observing the effect of the sub-methods designed for enhancing the performance of Algorithm
\ref{alg:circlePacking} one by one,  the efficiency of the algorithm
 used for initializing the upper bound is also explored if   the best-known value from
the literature are not fed into the algorithm. As a result,  the instances are solved by only using the upper bound
initialization algorithm from \cite{minRadius_Huang2006}, and these
results are compared with the results obtained by initializing the upper bounds as equal to the best-known
values in Table~\ref{tab:bestknowncomparison}. Algorithm ~\ref{alg:circlePacking}
terminates after 120 minutes of processing.  As seen in Table~\ref{tab:bestknowncomparison}, Algorithm~\ref{alg:circlePacking} performs good enough without knowing the best-known values by using the upper bound initialization algorithm stated in Section \ref{sec:otherimp}. Although it may not have found the best-known radii for some instances, the general performance of Algorithm~\ref{alg:circlePacking} is not effected conspicuously.

\begin{table}[H]
	\vspace{0.5cm}
		\caption{Effect of using best-known values within Algorithm~\ref{alg:circlePacking}.} \label{tab:bestknowncomparison} 
	\begin{center}
    \resizebox{\textwidth}{!}{
		\begin{tabular}{|c|cccc|cccc|}
     		\hline
     		& \multicolumn{4} {c|} {\bfseries Upper Bound Initialization Algorithm} & \multicolumn{4} {c|} {\bfseries Best-known Values}  \\ 
     		Instance	&	\shortstack{Upper\\Bound} &	\shortstack{Pre-processing \\Time (m)} &	\shortstack{Total Sol.\\Time (m)}	&	\shortstack{GAP (\%)} &	\shortstack{Upper\\Bound} &	\shortstack{Pre-processing \\Time (m)} &	\shortstack{Total Sol.\\Time (m)}	&	\shortstack{GAP (\%)} \\  \hline
			$Zimm-12$ & 28.371 & 3.72 & 14.43 & 0.59 & 28.371 & 2.71 & 12.92 & 0.57 \\
			$Zimm-13$ & 31.546 & 4.90 & 20.43 & 0.99 & 31.545 & 3.91 & 17.26 & 0.96 \\
			$Zimm-14$ & 35.096 & 6.04 & 37.05 & 0.43 & 35.096 & 5.41 & 36.43 & 0.37 \\
			$Zimm-15$ & 38.839 & 7.88 & 54.26 & 0.92 & 38.838 & 6.40 & 48.13 & 0.90 \\
			$Zimm-16$ & 42.457 & 9.94 & 63.56 & 0.85 & 42.457 & 8.51 & 61.30 & 0.84 \\
			$Zimm-17$ & 46.291 & 12.93 & 73.59 & 0.68 & 46.291 & 10.71 & 68.73 & 0.68 \\
			$Zimm-18$ & 50.129 & 15.07 & 99.26 & 0.83 & 50.120 & 12.87 & 90.36 & 0.46 \\
			$Zimm-19$ & 54.240 & 19.98 & 111.51 & 0.93 & 54.240 & 17.19 & 108.09 & 0.93 \\
			$Zimm-20$ & 58.401 & 24.91 & TLE & 1.97 & 58.401 & 21.63 & TLE & 1.97 \\
			\hline
		\end{tabular}
        }
	\end{center}
\end{table}

According to Table~\ref{tab:bestknowncomparison}, one can also say that Algorithm~\ref{alg:circlePacking} can solve problems up
to 19 circles within 120 minutes time limit whether or not the best-known is assumed to be given. In addition, the solution quality can be increased by solving the problem
containing  20 circles with a longer time limit. Since the present algorithm is a global optimization method, the solution quality
is guaranteed, especially for such a problem for which 
a tight lower bound is hard to verify. Therefore, it is possible to say that the proposed solution approach is quite successful as a global optimization method for CPP.

\subsection{Comparison with an Algorithm Designed for Equal Circle Case}\label{compSt:EqualCase}

In this subsection,   the results  for the equal circle instances are given.  
 For the packing of $n$-equal circles into a larger circle problem, a tailored algorithm is proposed by \cite{GlobalSolver_HUANG2011}, and the present approach is    compared against their algorithm. 
It can be  noted that   the proposed algorithm in \cite{GlobalSolver_HUANG2011} combines a local-search procedure and an improvement heuristic.  
The comparison of Algorithm~\ref{alg:circlePacking} with the algorithm introduced in \cite{GlobalSolver_HUANG2011} is given in Table~\ref{Tab:HuangOurAlgorithm}. The upper bounds are obtained by the upper bound initialization algorithm of \cite{minRadius_Huang2006}.
\begin{table}[H]
	\vspace{0.5cm}
		\caption{Results obtained by Algorithm~\ref{alg:circlePacking} and algorithm in \cite{GlobalSolver_HUANG2011}.} \label{Tab:HuangOurAlgorithm} 
	\begin{center}
    \resizebox{0.85\textwidth}{!}{
		\begin{tabular}{|c|ccc|cc|}
     		\hline
     		& \multicolumn{3}{c|}{\bfseries Algorithm~\ref{alg:circlePacking}}  & \multicolumn{2}{c|}{\bfseries Algorithm in \cite{GlobalSolver_HUANG2011}}\\ 
     		No.	& \shortstack{Upper Bound} &	\shortstack{Sol. Time (m)}	&	\shortstack{GAP (\%)}	&	\shortstack{Upper Bound} &	\shortstack{Sol. Time (m)} \\  \hline
			$Eq-20$ & 5.122 & 4.17 & 0.64 & 5.122 & 4.23 \\
			$Eq-25$ & 5.753 & 6.34 & 0.91 & 5.753 & 4.14 \\
			$Eq-30$ & 6.198 & 9.43 & 0.62 & 6.198 & 5.72 \\
			$Eq-35$ & 6.696 & 12.08 & 0.97 & 6.697 & 7.76 \\
			$Eq-40$ & 7.123 & 14.42 & 0.89 & 7.124 & 10.59 \\
			\hline
		\end{tabular}
        }
	\end{center}
\end{table} 
Algorithm~\ref{alg:circlePacking} presents an upper and a lower bound for the radius of the surrounding
circle. Hence, it is  ensured that the obtained circle is at most 1\% away from the
optimal-sized circle in which all circles are packed. However, the algorithm introduced
in Huang \& Ye (2011) only compares the obtained solutions with the best-known
values given in \cite{packomania}. Although this algorithm proposed by Huang \& Ye (2011)
solves the instances in shorter times, the increase in the solution time is not crucial. It is noted that their algorithm is tailored for equal circle case whereas the present algorithm is generic.

It can be also observed that Algorithm~\ref{alg:circlePacking} can solve the equal circle problem for the instances with a larger number of circles as well as the solution times are shorter than
the unequal circle problem. This observation is likely a result of the tighter
initial lower bounds. Since the minimum idle regions between the circles do not
change according to different configurations of circles for the equal circle case,  
the obtained lower bounds are tighter for the equal circle packing problem.

\subsection{Comparison with an Algorithm Tailored for a Rectangular Container}\label{compSt:UnequalCase}

In this subsection, the present algorithm is compared with another algorithm proposed by \cite{GlobalSolver_Stoyan2008} in which   the surrounding container is a rectangular strip. To handle this issue,   the feasible region definitions in Algorithm~\ref{alg:identifyingFeasibleRegions} are changed within the surrounding container and the definitions of the sets $\mathcal{L}_c, \mathcal{S}_c$.  

During the modification process, the guiding points are included by the rectangular strip whose width is known. 
 To ensure that a corresponding circle is totally included by the surrounding rectangle, guiding points are included by a rectangular ring whose outer boundary is $r_c$ units away from the boundary for circle $c$. In this problem, the aim is  to find the minimum length of the strip. In addition, the idle regions between the circles and the surrounding container should be updated accordingly. For that, the idle region is equal to the region between the tangent passing through both circles. The reduced regions are calculated according to the shape of the container. For the two largest circles, the first circle's center is located at the right upper part of the rectangle, and the second circle's center is located on the upper part of the diagonal connecting the right upper and left lower corner points of the rectangle. The set definitions within Algorithm~\ref{alg:identifyingFeasibleRegions} are accordingly changed. 
With these changes above,  Algorithm~\ref{alg:circlePacking} is  compared with the algorithm proposed in \cite{GlobalSolver_Stoyan2008}. During the comparisons, the initial upper bounds are obtained by the algorithm introduced in \cite{minRadius_Huang2006} which is also updated   to pack the circles into a rectangular strip.
\begin{table}[H]
	\vspace{0.5cm}
		\caption{Results obtained by Algorithm~\ref{alg:circlePacking} and the algorithm   in~\cite{GlobalSolver_Stoyan2008}.} \label{Tab:StoyanOurAlgorithm} 
	\begin{center}
    \resizebox{\textwidth}{!}{
		\begin{tabular}{|c|c|ccc|cc|}
     		\hline
     		& &	\multicolumn{3}{c|}{\bfseries Algorithm~\ref{alg:circlePacking}}  & \multicolumn{2}{c|}{\bfseries Algorithm in \cite{GlobalSolver_Stoyan2008}}\\ 
     		No.	& 	\shortstack{Strip \\Size (units)}  &	\shortstack{Upper\\ Bound} &	\shortstack{Solution\\ Time (m)}	&	\shortstack{GAP (\%)}	&	\shortstack{Upper \\ Bound} &	\shortstack{Solution \\ Time (m)}	\\  \hline
			$SY1$ & 9.5 & 17.291 & 47.35 & 0.91 & 17.461 & 34.14 \\ 
			$SY2$ & 8.5 & 14.375 & 23.07 & 0.37 & 15.604 & 8.11 \\ 
			$SY2-1$ & 9.0 & 13.572 & 24.16 & 0.12 & 13.653 & 7.64 \\ 
			$SY2-2$ & 9.5 & 12.758 & 19.83 & 0.76 & 12.547 & 6.01 \\ 
			$SY2-3$ & 11.0 & 11.163 & 20.49 & 0.85 & 11.214 & 6.69 \\ 
			$SY3$ & 9.0 & 14.321 & 31.29 & 0.82 & 15.171 & 18.09 \\ 
			$SY3-1$ & 8.5 & 15.904 & 36.40 & 0.77 & 16.463 & 17.26 \\ 
			$SY3-2$ & 9.5 & 13.691 & 34.05 & 0.94 & 13.713 & 13.52 \\ 
			$SY3-3$ & 11.0 & 11.564 & 29.52 & 0.37 & 11.827 & 16.52 \\ \hline
		\end{tabular}
        }
	\end{center}
\end{table} 
According to the results given in Table~\ref{Tab:StoyanOurAlgorithm}, the algorithm   in \cite{GlobalSolver_Stoyan2008} performs worse in terms of the solution quality, and the optimality
gap is unknown for these solutions. \cite{GlobalSolver_Stoyan2008} state that the proposed
algorithm   depends on the initialization procedure, and the initial solution affects
the solution quality of the algorithm. Also, their algorithm requires an
initial solution with all corresponding coordinates for the circles; however,  Algorithm~\ref{alg:circlePacking} requires an upper bound and a lower bound for the strip size during the
initialization procedure which can be assigned as the sum of the radii of the circles
and zero, respectively. Although the initial upper and lower bounds are not well-defined, the performance of 
Algorithm~\ref{alg:circlePacking} is not affected adversely since it is a bisection-type
algorithm and remove the half of the possible values at each iteration. However,
when comparing the algorithms in terms of the solution times, one can say that
the proposed procedure requires more time to solve the instances. This is a somewhat expected result since  the strength of
Algorithm~\ref{alg:circlePacking} is in guaranteeing the solution quality whereas \cite{GlobalSolver_Stoyan2008}
do not propose any optimality gap during the solution procedure. In addition, for many
of the instances, the present approach is able to find better upper bounds than the algorithm in
\cite{GlobalSolver_Stoyan2008}.

\section{Conclusion}\label{sec:conclusion}

In this article, the circle packing problem with the objective of minimizing the radius of a circular container is   studied. This nonconvex problem arises in many different settings and quite challenging to solve in general. Since global solvers perform poorly, a discretization-based solution approach is  designed to solve the problem to global optimality. The present solution approach iterates between a restricted and a relaxed version of the problem, which are both formulated as integer linear programming models, and terminate once the upper and lower bounds are within the user specified tolerance. This algorithm is enhanced using feasible region reduction, bound tightening and variable elimination strategies, and tested with a variety of instances from literature. The  computational experiments show that the proposed approach is orders-of-magnitude faster than the global solvers. Moreover, the proposed generic approach is competitive against heuristics tailored for specific cases, which cannot provide any optimality guarantee. As a future work,   real-life instances such as the 162-circle instance arising in the  automobile industry \citep{MinRadiusofContainer_Sugihara2004DiskPF} can be solved.

\subsection*{Acknowledgments} 
This work was supported by the  Scientific and Technological Research Council of Turkey under grant number 120M345.

\subsection*{Data Availability Statement} 
The input data of this study is collected from   the references \cite{packomania}, \cite{zimmcontest}, \cite{GlobalSolver_HUANG2011} and   \cite{GlobalSolver_Stoyan2008}.

\subsection*{Disclosure Statement} 
No potential conflict of interest was reported by the authors.

\bibliographystyle{tfcad}
\bibliography{references}

\end{document}